\documentclass[12pt]{amsart}
\usepackage{amsfonts,amsmath,amsthm,amscd,amssymb,latexsym,cite,verbatim,texdraw,floatflt,
caption2}
\RequirePackage{hyperref}

\usepackage[a4paper]{geometry}

\theoremstyle
{plain}

\begin{document}

\title{Decompositions of set-valued mappings }

\author{Igor Protasov}

\maketitle

\begin{flushright}
  {\it    On  100th anniversary of Professor V.S. $\check{C}$arin} \par
  \end{flushright}

\vskip 15pt

{\bf Abstract.}
Let $X$  be a  set,  $B_{X}$  denotes  the family of all subsets of    $X$   and
 $F: X \longrightarrow B_{X}$    be a set-valued mapping such that
   $x \in F(x)$,
 $sup  _{x\in X} |  F(x)|< \kappa$,  $sup  _{x\in X}  | F^{-1}(x)|< \kappa$
    for all  $x\in X$  and some  infinite cardinal $\kappa$.
Then  there exists a family $\mathcal{F}$  of  bijective  selectors of $F$  such
 that
 $|\mathcal{F}|<\kappa$ and
 $F(x) = \{ f(x): f\in\mathcal{F}\}$
   for each $x\in X$.
We apply this result to    $G$-space representations of balleans.

\vskip 10pt

\vskip 10pt


\vskip 10pt

Victor Sil'vestrovich  $\check{C}$arin is known as  the  founder  of topological algebra in Kyiv University,  but his mathematical  interests were not bounded by topological groups.  He encouraged and supported the activity of students and collaborators in many areas, in particular, in combinatorics.

\section{Decompositions}


For a set $X$, $B_{X}$  denotes the family of all subsets of $X$.
Given a set-valued
mapping
 $F: X\longrightarrow B _{X}$,  any
 function
  $f: X\longrightarrow X$   such that, for each $x\in X$,
   $f(x) \in  F(x)$  is  called a  {\it selector} of $F$.
We say that a selector $f$  is  {\it  bijective}  if $f: X\longrightarrow X$   is a bijection.
For $x\in X$,   we denote  $F ^{-1} (x) =\{ y\in X : x\in  F(y) \}$.

In section 1 we prove the mail result and apply   it to  $G$-space representations of  balleans  in section 2.

{\bf Theorem 1}. {\it
Let  $F :  X\longrightarrow  B _{X}$ be  a  set-valued  mapping  such  that
 $x\in F(x)$,   $sup  _{x\in X}  |F(x)| < \kappa$,
 $sup_{x\in X}  |F^{-1}(x)| < \kappa$
      for each  $x\in X$
      and  some  infinite cardinal $\kappa$.
Then there exists a family  $\mathcal{F}$   of bijective selectors of $X$   such  that
 $|\mathcal{F}| < \kappa $  and
 $F (x) =\{ f(x) : f\in \mathcal{F} \}$
    for  each $x\in X$.

\vspace{5 mm}

Proof.} We consider two cases.
\vspace{5 mm}

{\it Case  $\kappa = \omega$}.
We put $\mathcal{P} = \{ F(x):  x\in X\}$  and  define  a graph $\Gamma$  with the  set of vertices  $\mathcal{P}$   and    the  set of   edges
  $\{\{ F(x),   F(y) \} : F(x)\cap F(y) \neq\emptyset\}$.
We  take a natural  number $m$   such  that
 $m>  sup _{x\in X}|F(x)|$,
 $m> sup |F^{-1}(x)|$ and show that the  local  degree  of each vertices of
 $\Gamma$
  does not exceed  $m ^{2} - 1$.
Assume the contrary and choose  $y\in X$ and distinct
$y_{1},\ldots , y_{m ^{2}}  \in X$ such  that
$F(y)\cap F(y_{i}) \neq\emptyset\}$
  for every  $i\in  \{1,\ldots , m^{2} \}$.
Then $y_{i} \in F ^{-1}  F(y)$  but, by  the choice of $m$, we  have
$|F ^{-1}  F(y)|< m^{2}$.

We use the following simple fact  \cite{b2}:  if the  local degree of each vertices of a graph
$\Gamma^{\prime}$
does not exceed $k$ then the chromatic number of $\Gamma^{\prime}$  does not exceed $k + 1$.

Hence the set  $\mathcal{P}$ of vertices of $\Gamma$
 can be partition
 $\mathcal{P}_{1},\ldots  , \mathcal{P}_{m^{2}}$
   so that any two vertices  from  each
   $\mathcal{P}_{i}$
   are  not incident.

To construct the family $\mathcal{F}$,  we enumerate
$\mathcal{P}_{i}= \{F(y_{\alpha}): \alpha< \gamma  \}.$
Let $M= sup _{x\in X}|F(x)|$.
Then we enumerate each
$F(y _{\alpha})$  (with repetitions,  if necessary)
$F(y _{\alpha})= \{y_{\alpha j}): j< M\},$
$y_{\alpha_{0}}  =  y_{\alpha} $.
For each
$j< M$,  we define a  bijective  function $f _{j}$  such that $f _{j}$
 acts as a transposition of $y_{\alpha}$  and
 $y_{\alpha j}$
  at each $F(y_{\alpha})$  and identically at all other   elements of $X$.
We put
$\mathcal{F}_{i}=  \{ f_{j}: j< M \}$
  and note that
  $\mathcal{F}= \mathcal{F}_{1} \cup  \ldots  \cup \mathcal{F}_{m^{2}} $
   is the desired family of  selectors of $F$.
\vspace{5 mm}

{\it Case  $\kappa > \omega$}.
We take an  infinite  cardinal
$\sigma$
 such that $\sigma < \kappa$  and
 $|F(x)|\leq\sigma$,  $|F^{-1}(x)|\leq\sigma$
  for each  $x\in X$.
Then we define  a partition $\mathcal{P}$ of  $X$  such that  each  $P\in \mathcal{P}$
is  the minimal by inclusion subset of $X$  satisfying
$F(y)\in P$, $F^{-1}(y)\in P$
 for each $y\in P$.
Constructively, every $P$  can be obtained applying   to  $x\in P$  the sequence of operations
$F$, $F^{-1}:  F(x)$, $F^{-1} F(x)$,  $FF^{-1}  F(x), \ldots$.
Then $P$  is the union of all numbers of this sequence.

By the choice of $\sigma$, we have $|P|\leq \sigma$.
We enumerate  $\mathcal{P}=\{ P_{\alpha}: \alpha< \gamma \} $,
$P_{\alpha}=\{ x_{\alpha j}: j< \gamma \}.$
For each $j< \sigma$,  we choose a family
$\mathcal{F}_{j}$
 of bijective selectors of $F$ such that
  $|F_{j}|\leq\sigma$  and
  $F(x_{\alpha j})= \{ f(x_{\alpha j}): f\in \mathcal{F}_{j} \}$
    for each
    $\alpha<\gamma$,
      see the case $\kappa=\omega$.
Then
$\bigcup_{j< \sigma} \mathcal{F}_{j}$
is the desired family  $\mathcal{F}$   of  bijective selectors   of  $F$.
$ \  \  \Box$

\section{Applications}

Let $X$  be a set.  A family $\mathcal{E}$  of subsets of $X\times X$ is called a {\it coarse structure}  if
\vspace{5 mm}

\begin{itemize}
\item{}   each $E\in \mathcal{E}$  contains the diagonal  $\bigtriangleup _{X}$,
$\bigtriangleup _{X}= \{(x,x): x\in X\}$;
\vskip 5pt

\item{}  if  $E$, $E^{\prime} \in \mathcal{E}$ then $E\circ E^{\prime}\in\mathcal{E}$ and
$E^{-1}\in \mathcal{E}$,   where    $E\circ E^{\prime}=\{(x,y): \exists z((x,z) \in  E,  \   \ (z, y)\in E^{\prime})\}$,   $E^{-1}=\{(y,x): (x,y)\in E\}$;
\vskip 5pt

\item{} if $E\in\mathcal{E}$ and $\bigtriangleup_{X}\subseteq E^{\prime}\subseteq E  $   then
$E^{\prime}\in \mathcal{E}$;
\vskip 5pt

\item{}  for any   $x,y\in X$, there exists $E\in \mathcal{E}$   such that $(x,y)\in E$.

\end{itemize}
\vskip 7pt

A subset $\mathcal{E}^{\prime} \subseteq \mathcal{E}$  is called a
{\it base} for $\mathcal{E}$  if, for every $E\in \mathcal{E}$, there exists
  $E^{\prime}\in \mathcal{E}^{\prime}$  such  that
  $E\subseteq E ^{\prime}$.
For $x\in X$,  $A\subseteq  X$
we denote $E[x]= \{y\in X: (x,y) \in E\}$, $E[A]= \cup_{a\in A} E[a]$
 and say $E[x]$ and $E[A]$
 are {\it balls of radius $E$
   around} $x$  and $A$.

The pair $(X,\mathcal{E})$ is called a {\it coarse space} \cite{b6}  or a {\it ballean} \cite{b5}.

Let $(X,\mathcal{E})$, $(X^{\prime},\mathcal{E}^{\prime})$ be coarse spaces.
A mapping $f: X\longrightarrow X^{\prime}$ is called {\it macro-uniform} if, for every
$E\in \mathcal{E}$ there exists $E^{\prime}\in \mathcal{E}^{\prime}$
such that
$E[x]\subseteq E^{\prime}[f(x)]$.
If $f$ is a bijection such that $f, f ^{-1}$  are macro-uniform then $f$ is called an {\it asymorphism}.

Now we describe some general way of  constructing balleans.
Let $G$ be a group.
A family $\mathcal{I}$ of subsets of $G$  is called an {\it ideal}  if, for every
 $A, B \in \mathcal{I}$ and $A^{\prime}\subseteq  A$,
  we have $A\cup B \in \mathcal{I}$ and $A^{\prime} \in \mathcal{I}$.
An ideal $\mathcal{I}$ is called a {\it group ideal}  if
$F\in \mathcal{I}$ for every finite subset of $G$ and $A, B\in \mathcal{I}$  imply $AB^{-1} \in \mathcal{I}$.

Let a group $G$ acts transitively on a set $X$  by the rule
$(g, x)\longmapsto  gx$, $g\in X$, $x\in X$.
Every group ideal $\mathcal{I}$ on $G$  defines the ballean $(X, G, \mathcal{I})$  on $X$
with  the base of entourages
$\{\{ (x, y):  y\in Ax\}: A\in\mathcal{I}\}$.
By Theorem 1 from
\cite{b3}, for every ballean  $(X, \mathcal{E})$, there
 exist a group  $G$ of permutations of $X$  and  a group ideal  $\mathcal{I}$  on $G$  such that  $(X, \mathcal{E})$  is asymorphic  to $(X, G, \mathcal{I})$.

\vspace{7 mm}

{\bf Theorem 2}. {\it
Let  $(X, \mathcal{E})$  be a ballean and let $\kappa$ be an  infinite cardinal such that,  for  each
$E\in \mathcal{E}$,  $sup_{x\in E} | E[x] |< \kappa$.
Then there exist   a group  $G$ of  permutations
of $X$ and a group ideal  $\mathcal{I}$  on $G$ such that  $(X, \mathcal{E})$ is asymorphic to  $(X, \mathcal{E}, \mathcal{I})$   and $|A|< \kappa$  for each  $A\in \mathcal{I}$.

\vspace{5 mm}

Proof.}
For each $E\in \mathcal{E}$,  we define a mapping $F _{E} : X\longrightarrow  B _{X}$  by $F _{E} (x)=E[x]$.
By Theorem 1, there exists a family $F _{E}$  of permutations of $X$  such  that $|\mathcal{F} _{E}|<\kappa $
and $F _{E} (x)=\{ f(x): f\in \mathcal{F}_{E}\}$  for each $x\in X$.
We denote by $\mathcal{I}$  the minimal by inclusion group ideal of $G$ such that
$\mathcal{F} _{E} \in \mathcal{I}$  for each $E\in \mathcal{E}$.
Then $(X, \mathcal{E})$  is asymorphic  to $(X, G, \mathcal{I})$.  $ \ \ \ \Box$

\vspace{5 mm}

In the case   $\kappa=\omega $,
Theorem 2 was proved in \cite{b4}.
For its applications see Remark 3.5  in \cite{b1}.

A ballean $(X, \mathcal{E})$ is called  {\it cellular}  if $\mathcal{E}$ has a base  consisting of equivalence  relations.
By Theorem 3 from \cite{b3}, every  cellular ballean  is asymorphic to some ballean  $(X, G, \mathcal{I})$  such that  $\mathcal{I}$ has a base consisting of  subgroups of $G$.

A ballean $(X, \mathcal{E})$  is called {\it finitary} if, for every $E\in \mathcal{E}$  there exists  a natural number $m$  such $|E[x]|< m $ for each $x\in X$.
The  finitary ballean of a $G$ space $X$  is the ballean  $(X, G, \mathcal{I})$, where $\mathcal{I}$  is the ideal of all finite subsets of $G$.
\vspace{7 mm}

{\bf Theorem 3}. {\it For every finitary cellular ballean  $(X, \mathcal{E})$ there exists  a locally finite group of permutations of $X$  such that $(X, \mathcal{E})$  is asymorphic  to the  finitary  ballean of  $G$-space  $X$.

\vspace{5 mm}

Proof.}
We take a base $\mathcal{E}^{\prime}$ of consisting  of  partitions  of $X$.
 For every $\mathcal{P}\in \mathcal{E}$  we  pick a natural number  $n _{\mathcal{P}}$ such  that $|P|\leq n _{\mathcal{P}} $  for each $P\in \mathcal{P}$.
We denote by   $G _{\mathcal{P}}$  the direct product of the family of symmetric groups $\{S_{m}: m\leq n _{\mathcal{P}}\}$   and note that
$G _{\mathcal{P}}$
 acts on each $P\in \mathcal{P}$  so that $G _{\mathcal{P}} x = P$  for  each $x\in P$.
Then the group  $G$  generated  by the family
$\{G _{\mathcal{P}}: \mathcal{P}\in\mathcal{E}^{\prime}\}$
   satisfies the conclusion of Theorem 3. $ \ \ \ \Box$

\vspace{6 mm}

CONTACT INFORMATION

I.~Protasov: \\
Faculty of Computer Science and Cybernetics  \\
        Kyiv University  \\
         Academic Glushkov pr. 4d  \\
         03680 Kyiv, Ukraine \\ i.v.protasov@gmail.com

\end{document}